\documentclass{article}
\usepackage{amsthm}
\usepackage{bbm}
\usepackage[title,titletoc]{appendix}
\usepackage{lineno}
\usepackage{graphicx}
\usepackage{amsmath,amssymb,amsthm,amsfonts}
\usepackage{amssymb}
\usepackage{color}
\usepackage{amsfonts}
\usepackage{amssymb}
\usepackage{amsthm}
\usepackage{amsmath}
\usepackage{empheq}
\usepackage{indentfirst}
\usepackage{cite}
\usepackage{mathrsfs}
\usepackage{lineno}
\allowdisplaybreaks[4]
\newtheorem*{thmA}{Theorem A}
\newtheorem*{thmB}{Theorem B}
\newtheorem{thm}{Theorem}[section]
\newtheorem{cor}{Corollary}[section]
\newtheorem{lem}{Lemma}[section]
\newtheorem{prop}{Proposition}[section]

\theoremstyle{definition}

\newtheorem{rem}{Remark}[section]

\numberwithin{equation}{section}

\DeclareMathSymbol{\C}{\mathalpha}{AMSb}{"43} \topmargin-.1in

\def\O{\mathcal{O}}

\def\Om{\Omega}
\def\r{\mathbb{R}}
\def\f{\frac}
\def\F{\displaystyle\frac}
\def\dis{\displaystyle}
\def\Io{\dis\int_{\Omega}}
\def\io{\int_{\Omega}}
\def\D{\Delta}
\def\e{\epsilon}
\newcommand{\beq}{\begin{equation}}
\newcommand{\eeq}{\end{equation}}
\newcommand{\bea}{\begin{array}}
\newcommand{\eea}{\end{array}}

\newcommand{\lam}{\lambda}
\newcommand{\Lam}{\Lambda}
\newcommand{\alp}{\alpha}

\newcommand{\R}{\mathbb{R}}
\newcommand{\gam}{\gamma}
\newcommand{\Gam}{\Gamma}
\newcommand{\del}{\delta}

\newcommand{\al}{\alpha}
\newcommand{\be}{\beta}
\def\th{\theta}
\newcommand{\pa}{\partial}
\def\dx{\mathrm{d}x}
\def\dt{\mathrm{d}t}
\def\ds{\mathrm{d}s}

\def\ni{\noindent}
\def\proof{{\ni\bf Proof:\quad}}
\def\proofend{{\hfill$\Box$}\\}
\def\prooftd1{{\ni\bf Proof of Theorem \ref{large initial data}: }}
\def\proofprop2{{\ni\bf Proof of Theorem \ref{thm1}: }}

\def\proofconrate{{\ni\bf Proof of Theorem \ref{Convergence rate}: }}

\def\proofthA{{\ni\bf Proof of Theorem A: }}
\def\proofthB{{\ni\bf Proof of Theorem B: }}

\title{Asymptotic and Quenching Behavior for a family of Parabolic System with General Singular Nonlinearities}
\author{Qi Wang
\thanks{College of Science, University of Shanghai for Science and Technology,
Shanghai 200093, P.R. China.
Email: \texttt{qwang@usst.edu.cn}.
},
Yanyan Zhang
\thanks{Corresponding author. School of Mathematical Sciences,
Shanghai Key Laboratory of Pure Mathematics and Mathematical
Practice, East China Normal University, Shanghai 200241, P.R. China.
Email: \texttt{yyzhang@math.ecnu.edu.cn}.
Y. Y. Zhang is sponsored
by ``Chenguang Program" supported by Shanghai Educational
Development Foundation and Shanghai Municipal Education Commission
[grant number: 13CG20]; NSFC [grant number: 11431005]; and STCSM
[grant number: 18dz2271000].} }
\date{\today}

\begin{document}

\maketitle

\begin{abstract}
This study is concerned with  a family of parabolic system with general
singular nonlinearities, which is a generalization of MEMS system.
To some extent, the classification of global existence and
  quenching according to parameters and initial data is
given. Moreover, the convergence rate is also obtained. We
point out that compared to single MEMS equation, new ideas and
techniques are introduced in obtaining the convergence rate for system
in our study. In fact, due to the lack of variational
characterization for the first eigenvalue of the linearized elliptic
system, the methods in obtaining convergence rate for single
equation cannot work completely here.
\end{abstract}



\section{Introduction}

This study is concerned with the following  coupled generalized singular
parabolic system of the form \beq \tag{P} \label{P}
\left\{\arraycolsep=1.5pt
\begin{array}{lll}
u_t-\Delta u=\displaystyle\lambda\alpha(x)f(v),\   \ &\hbox{in} \,\  \Omega\times(0,T),\\[2mm]
v_t-\Delta v=\displaystyle\mu\beta(x)g(u), \   \ &\hbox{in} \,\  \Omega\times(0,T),\\[2mm]
u=v=0,\ &\hbox{on} \ \partial \Omega\times(0,T),\\[2mm]
u(x,0)=u_0(x),\ v(x,0)=v_0(x), \ &\hbox{in}\,\  \Omega,
\end{array}\right.
\eeq where $\Omega\subset\mathbb{R}^N$ is a smooth bounded domain,
$\lambda$ and $\mu$ are positive parameters, $\alpha(x)$ and
$\beta(x)$ are nonnegative nontrivial H\"{o}lder continuous
functions in $\bar{\Omega}$, $f,g$ satisfy \beq \label{hhhhh}
\tag{H1}
\begin{array}{l}
f,g\in C^{2}[0,1)~\mbox{are positive, increasing and strictly convex such that}~\lim\limits_{v\rightarrow 1^-}f(v)=\lim\limits_{u\rightarrow 1^-}g(u)=+\infty,\\
\end{array}
\eeq and the initial data satisfy \beq \label{assumption on initial
data}\tag{H2} u_{0}(x),v_{0}(x)\in C^2(\bar{\Om}),0\le u_0,v_0<1,~
u_0=v_0=0~\mbox{on}~\partial\Omega. \eeq In fact $f(\cdot),g(\cdot)$
can be one of the following nonlinearities: $1-\ln(1-\cdot)$,
$e^{\frac{1}{(1-\cdot)}}$, $\frac{1}{(1-\cdot)^p}$ with $p>0$.
\begin{rem}
In \eqref{hhhhh}, we  fix the blow up level at $u=1, v=1$ for
simplicity.
 It is easy to see that with the scaling, our approaches work for $f,g$ blowing up at any positive values $a$ and $b$, respectively.
\end{rem}
Recall that the scalar equation \beq\label{single}
\left\{\arraycolsep=1.5pt
\begin{array}{lll}
u_t-\Delta u=\displaystyle\lambda\alpha(x)f(u),\   \ &\hbox{in} \,\  \Omega\times(0,T),\\[2mm]
u=0,\ &\hbox{on} \ \partial \Omega\times(0,T),\\[2mm]
u(x,0)=u_0(x)\in[0,1) \ &\hbox{in}\,\  \Omega
\end{array}\right.
\eeq  as well as the associated stationary equation \beq
\label{S-single} -\Delta u=\displaystyle\lambda\alpha(x)f(u)\
\hbox{in} \,\  \Omega,\ u=0\ \hbox{on} \ \partial \Omega \eeq with
$f$ satisfying \eqref{hhhhh}  have been studied in \cite{YZ}. More
precisely, it is  showed  in \cite{YZ} that for any given
$\alpha\geq0$ and  $f$ satisfying \eqref{hhhhh}, there exists a
critical value $\lambda^*>0$ such that if $\lambda>\lambda^*$, no solution of \eqref{S-single} exists, and the
solution to \eqref{single} will reach the value 1 at finite time
$T$, i.e., the so called quenching or touchdown phenomenon occurs;
while for
$\lambda\in(0,\lambda^*)$,
problem  \eqref{S-single} is solvable and the solution to
\eqref{single} is global and convergent for some $u_0$.
Moreover, the convergence rate is obtained in \cite{LL}, by noting that the first eigenvalue of the linearized elliptic equation having a variational characterization.

In fact, besides \cite{YZ,LL}, for the particular case
$f(u)=(1-u)^{-p},p>0$, especially for $p=2$, as the mathematical
model of micro-electromechanical systems (MEMS), \eqref{single} has
been extensively studied by many authors in recent
years(cf.\cite{Esposito,GZZ,GG} and references therein). MEMS device
consists of an elastic membrane suspended over a rigid ground plate.
For MEMS, $u$ denotes  the normalized distance between the membrane
and the ground plate,  $\alpha(x)$ represents the permitivity
profile. When a voltage $\lambda$ is applied, the membrane deflects
toward the ground plate and a snap-through may occur when it exceeds
a certain critical value $\lambda^*$ (pull-in voltage). This creates
a so-called “pull-in instability”, which greatly affects the
design of many devices (cf.\cite{Pelesko,Esposito} for more
details).

As for system \eqref{P},  if
$f(v)=(1-v)^{-p},g(u)=(1-u)^{-q},p,q>0,$ due to the reason above,
system \eqref{P} is called  general MEMS system(see \cite{Fazly}).
For system \eqref{P},  some sufficient conditions  for global
existence,  quenching and quenching time estimates of
solutions, as well as the non-simultaneous quenching criteria are
studied in
\cite{BNN,Jia,Gu,Selcuk,zhengsining,mu,pei} and the  references cited
therein. However, as far as we know,  the classification of global
existence and quenching is not fully described, and the
convergence rate for global solution is not considered before.

In this paper, to some extent, we  show the  classification  of  global  existence and  quenching  according  to  parameters  and  initial  data. Then we further study the global solution's asymptotic behavior, such as convergence and convergence rate.
The
solution $(u,v)$ of \eqref{P} is called quenching at time
$t=T<+\infty$ if \beq \label{defn of quenching}
\limsup\limits_{t\rightarrow T^-}
(\max\{\max\limits_{\Om}u(\cdot,t),\max\limits_{\Om}v(\cdot,t)\})=1.
\eeq
Similarly as the considerations in single MEMS equation (cf. \cite{YZ,Esposito}), there exists a close relationship between
system \eqref{P} and
the associated stationary problem \beq \tag{E} \label{E}
\left\{\arraycolsep=1.5pt
\begin {array}{lll}
-\Delta w=\displaystyle\lambda\alpha(x)f(z),\ \ &\hbox{in} \,\  \Omega,\\[2mm]
-\Delta z=\displaystyle\mu\beta(x)g(w), \ \ &\hbox{in} \,\  \Omega,\\[2mm]
w=z=0,\ &\hbox{on} \ \partial \Omega.
\end{array}\right.
\eeq
Therefore,
we first consider the associated stationary problem
\eqref{E}.

\begin{figure}[htbp]
\centering
{\scalebox{0.5}[0.5]{\includegraphics{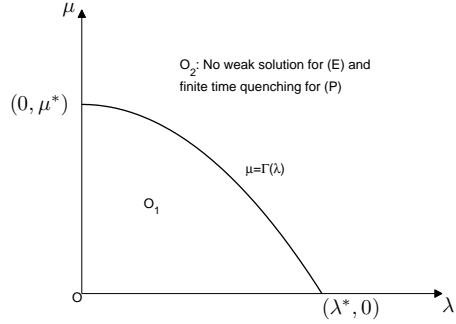}}}
\caption{\label{fig:shiyitu}
\em The critical curve $\Gamma$ in $(\lam,\mu)$-plane}
\end{figure}

Recall that for \eqref{E} with $f(\cdot)=g(\cdot)=(1-\cdot)^{-2}$ ,
it has been proved in \cite{Clemente} that there exists a critical
curve $\Gamma$  splitting the positive quadrant of the
$(\lambda,\mu)$-plane into two disjoint sets $\O_1$ and $\O_2$ such
that the elliptic problem has a smooth minimal stable solution
$(w_{\lambda,\mu},z_{\lambda,\mu})$ for $(\lambda,\mu)\in\O_1$,
while for $(\lambda,\mu)\in\O_2$ there is no solution of any kind.
In this paper,  we note that these results can be extended
to general elliptic problem \eqref{E} as shown in Theorem
A, which can be  illustrated by Figure \ref{fig:shiyitu}.

\begin{thmA}
There exist $0<\lambda^*,\mu^*<+\infty$, and  a non-increasing
continuous curve $\mu=\Gamma(\lambda)$ connecting $(0,\mu^*)$ and
$(\lambda^*,0)$ such that  the positive quadrant $\r^+\times\r^+$ of
the $(\lam,\mu)$-plane is  separated into two connected components
$\O_1$ and $\O_2$. For $(\lam,\mu)\in\O_1$, problem \eqref{E} has a
positive classical minimal solution $(w_{\lam,\mu},z_{\lam,\mu})$.
Otherwise, for $(\lam,\mu)\in\O_2$, \eqref{E} admits no weak
solution.
\end{thmA}

Since Theorem A can be established similarly to
\cite{Clemente}, we sketch  the proof in Appendix \ref{Stationary
problem} for simplicity.

The second fundamental problem  is the local existence and
uniqueness of solution to \eqref{P}. In fact, due to the increasing monotonicity of  functions $f,g$, one can get the
following Theorem B via comparison principle, see such as
\cite[Theorem 13, Chapter3]{PW}
and  C. V. Pao\cite[Theorem 2.2]{Pao2}. For the convenience of
readers, we sketch  the proof in Appendix \ref{proof of local
existecnce and uniqueness}.
\begin{thmB}[Local existence and uniqueness]
\label{Local existence} Suppose \eqref{hhhhh}
 and \eqref{assumption on initial data} hold. Then for any $\lam>0,\mu>0$ there exists $T>0$ such that problem
\eqref{P} has a unique solution $(u,v)$ on $[0,T)$. Moreover $(u,v)\in
C^{1}((0,T),C^2(\Om,\mathbb{R}^2))$.
\end{thmB}

The main results of this paper are stated as follows. By the way, in this paper, we use $(u,v)\leq(w,z)$ to denote $u\leq w$ and $v\leq z.$ $\|\cdot\|_p$ denotes  the standard norm of
$L^{p}(\Om)$.

\begin{thm}[Global existence vs quenching]
\label{Global existence vs quenching} Suppose \eqref{hhhhh}
and \eqref{assumption on initial data}
  hold.
 Let $\mathcal{O}_1$ and $\mathcal{O}_2$ be the connected component defined
in Theorem A as well as $(w_{\lam,\mu},z_{\lam,\mu})$ be the minimal
solution of \eqref{E}.

$(a)$ If $(\lam,\mu)\in\mathcal {O}_1$, then there hold:

\ \ \ \ $(a_1)$  For $(u_{0},v_{0})\le
(w_{\lam,\mu}, z_{\lam,\mu})$,  the unique solution
$(u(x,t),v(x,t))$ to \eqref{P} exists globally and converges uniformly
to $(w_{\lam,\mu},z_{\lam,\mu})$ as $t\rightarrow +\infty.$ If $(u_{0}(x),v_{0}(x))$ is further a subsolution of \eqref{E},  the convergence is monotone increasing.

\ \ \ \  $(a_{2})$ Suppose that \eqref{E} has solutions more than one. Let $(w_1,z_1)$ be any solution different from $(w_{\lam,\mu}, z_{\lam,\mu}).$ Then

\ \ \ \ \ \ \ \ $(a_{21})$ For $(u_{0},v_{0})\le\not\equiv
(w_1, z_1)$, the unique solution
$(u(x,t),v(x,t))$ to \eqref{P}  exists globally and converges  uniformly to
$(w_{\lam,\mu},z_{\lam,\mu})$ as $t\rightarrow +\infty.$ If there exists $s\in(0,1)$ such that $(u_{0},v_{0})=s(w_{\lam,\mu}, z_{\lam,\mu})+(1-s)(w_1, z_1),$
 the convergence is monotone decreasing.

\ \ \ \ \ \ \ \ \ $(a_{22})$ For $(u_{0},v_{0})\ge \not\equiv
(w_1, z_1)$,   the solution $(u,v)$ to \eqref{P} will quench (at a finite time or an infinite time).

$(b)$ If $(\lam,\mu)\in\mathcal {O}_2$,
  then  for any $(u_0,v_0)$ satisfying \eqref{assumption on initial data}, the solution $(u,v)$ to \eqref{P} will quench at a finite time.

$(c)$ For any given $\lam>0$ and $\mu>0$, once
$\dis\int_{\Om}u_{0}\phi\dx>\bigg({\lam_{1}\dis\int_{\Om}{\phi}/{\al}\dx}\bigg)/\big({\lam f(0)}\big)$ or $\dis\int_{\Om}v_{0}\phi\dx>\bigg({\lam_{1}\dis\int_{\Om}{\phi}/{\be}\dx}\bigg)/
 \big({\mu g(0)}\big)$,
 the  solution $(u,v)$ of \eqref{P}
  should quench at a finite time.  Here $\lam_{1}>0$ is the first eigenvalue of $-\Delta$ in $H_{0}^{1}(\Om)$   and $\phi$ is the corresponding eigenfunction satisfying $\dis\int_{\Om}\phi\dx=1$.
\end{thm}

\begin{rem}
\label{C1 converge}
If $1\leq n\leq 3,$ then all the convergence in Theorem \ref{Global existence vs quenching} further hold in $C^1$ norm.
\end{rem}

\begin{rem}
As far as we know, the multiplicity of  solution to \eqref{E} is  open.
\end{rem}

\begin{rem}
When the minimal solution $(w_{\lam,\mu},z_{\lam,\mu})$ is the only solution to \eqref{E}, except for the case in Theorem \ref{Global existence vs quenching} $(c)$, the behavior of solution  $(u(x,t),v(x,t))$ to \eqref{P} with initial data above $(w_{\lam,\mu},z_{\lam,\mu})$ depends on circumstances, which are complicated.
\end{rem}

\begin{cor} \label{estimates}
$(1)$ There   hold$$\int_{\Om}w_{\lam,\mu}\phi\dx\le\bigg({\lam_{1}\dis\int_{\Om}{\phi}/{\al}\dx}\bigg)/\big({\lam f(0)}\big)\quad \text{and}\quad  \int_{\Om}z_{\lam,\mu}\phi\dx\le\bigg({\lam_{1}\dis\int_{\Om}{\phi}/{\be}\dx}\bigg)/
 \big({\mu g(0)}\big).$$

$(2)$ If \eqref{E} has another solution $(w_1,z_1),$ then there further  hold
$$\int_{\Om}w_1\phi\dx\le\bigg({\lam_{1}\dis\int_{\Om}{\phi}/{\al}\dx}\bigg)/\big({\lam f(0)}\big) \quad \text{and}\quad  \int_{\Om}z_1\phi\dx\le\bigg({\lam_{1}\dis\int_{\Om}{\phi}/{\be}\dx}\bigg)/
 \big({\mu g(0)}\big).$$
\end{cor}

\begin{thm}[Convergence rate]
\label{Convergence rate} Suppose \eqref{hhhhh}
and \eqref{assumption on initial data}
  hold.  For cases in Theorem \ref{Global existence vs quenching} $(a_1)$ and $(a_{21})$, if
 $1\le n\le3$,
 the convergence rate can be estimated by
 \beq\label{rate}
  \|u(x,t)-w_{\lam,\mu}(x)\|^2_2+\|v(x,t)-z_{\lam,\mu}(x)\|^2_2\leq C_0 \exp\bigg({-\min\bigg\{2\lambda_1,{\nu_1}/{2}\bigg\}t}\bigg),\quad \ t>T_0
  \eeq
  for some $T_0>0$, where $C_0$ is a constant depending at most on $u_0,v_0,w_{\lam,\mu},z_{\lam,\mu},w_1,z_1,\varphi_1,\psi_1.$
Here
$\lambda_1>0$ is the first eigenvalue of $-\Delta$ on
$H_0^1(\Omega)$,
 $\nu_1>0$ is the first eigenvalue of linearized elliptic
system \eqref{7}, and $(\varphi_1,\psi_1)$ are the corresponding
strictly positive eigenfunction defined in Lemma \ref{aaa6}.
  \end{thm}

 We  remark that as one of  main contribution of this paper,
obtaining the  convergence rate \eqref{rate} needs new ideas and
techniques. 
In fact, for single
parabolic MEMS equation, in obtaining the convergence rate of global solution \cite{LL}, the first eigenvalue of the linearized elliptic equation having a variational characterization plays an important role. However,  no such analogous formulation is available for coupled system \eqref{P} considered in this paper (see \cite{Clemente}).
Another commonly used idea is $\L$ojasiewicz--Simon method.  Though there exists a
Lyapunov function for  \eqref{P}, however, it is not  coercive. Hence we can not apply $\L$ojasiewicz--Simon method  directly for system \eqref{P}, either.

This paper is organized as follows. In Section 2, we will show some preliminary results. Then Theorem \ref{Global existence vs quenching}-\ref{Convergence rate} will be proved in
Section  \ref{section Global existence vs quenching} and \ref{Convergence rate for MEMS system}, respectively. At last, we show the proof of Theorem A and B in Appendix A and B,
respectively.

\section{Preliminary results}\label{3.1}

The proof of Theorem \ref{Global existence vs quenching} relies on some results, which are presented as follows. Proposition \ref{triple} is concerned with the structure of the stationary solution set. Proposition \ref{ut vt} is concerned with the monotonicity of global solution to \eqref{P}.

\begin{prop} \label{triple}
There cannot be a triple $(w_i,z_i)(i=1,2,3)$ of solutions to  \eqref{E} with $(w_1,z_1)\ll(w_2,z_2)\ll(w_3,z_3)$. Here $w_1\ll w_2$ means that $\gamma \rho\leq w_2-w_1(x\in \bar{\Omega})$ for some positive number $\gamma$, $\rho(x)$ being the distance from $x$ to $\partial \Omega.$
\end{prop}
\proof
Suppose by contradiction that there exists a triple $(w_i,z_i)(i=1,2,3)$ of solutions to \eqref{E}, i.e.,
\beq \left\{\arraycolsep=1.5pt
\begin {array}{lll}
-\Delta w_i=\lam\al(x)f(z_i),\   \ &\hbox{ in} \,\  \Omega,\\[2mm]
-\Delta z_i=\mu\be(x)g(w_i), \   \ &\hbox{ in} \,\  \Omega,\\[2mm]
w_i=z_i=0,\ &\hbox{ on} \ \partial \Omega.
\end{array}\right.
\eeq
It then follows that
\beq
\label{3-2-1}
\left\{\arraycolsep=1.5pt
\begin {array}{lll}
-\Delta (w_3-w_2)&=\lam\al(x)(f(z_3)-f(z_2)),\   \ &\hbox{ in} \,\  \Omega,\\[2mm]
-\Delta (z_3-z_2)&=\mu\be(x)(g(w_3)-g(w_2)), \   \ &\hbox{ in} \,\  \Omega,\\[2mm]
-\Delta (w_2-w_1)&=\lam\al(x)(f(z_2)-f(z_1)),\   \ &\hbox{ in} \,\  \Omega,\\[2mm]
-\Delta (z_2-z_1)&=\mu\be(x)(g(w_2)-g(w_1)), \   \ &\hbox{ in} \,\  \Omega.
\end{array}\right.
\eeq
If we multiply the first equation in \eqref{3-2-1} by $z_2-z_1$, the second equation in \eqref{3-2-1} by $w_2-w_1$ and integrating over $\Om$,
then by \eqref{3-2-1} there exists $\xi_w\in(w_2,w_3),$ $\xi_z\in(z_2,z_3),$ $\eta_w\in(w_1,w_2,)$ $\eta_z\in(z_1,z_2)$ such that
\beq
\Io\mu\be(x)(w_3-w_2)(w_2-w_1)(g'(\eta_w)-g'(\xi_w))\dx+\Io\lam\al(x)(z_3-z_2)(z_2-z_1)(f'(\eta_z)-f'(\xi_z))\dx=0.
\eeq
This is impossible,
since by the strict convexity of $f,g$ there hold $g'(\eta_w)-g'(\xi_w)<0$ and $f'(\eta_z)-f'(\xi_z)<0$  by noting that $\eta_w<\xi_w, \eta_z< \xi_z$ in $\Omega.$ The proof is completed.
\proofend

\begin{rem}
Suppose that $(w_i,z_i)(i=1,2)$ are solutions to \eqref{E}. If $(w_1,z_1)\leq \not\equiv (w_2,z_2),$ then we have $(w_1,z_1)\ll (w_2,z_2)$ by Hopf Lemma.
\end{rem}

To
verify Proposition \ref{ut vt} below, we need to borrow a comparison
principle for the parabolic system below, which can be derived from
\cite[Theorem 13, Chapter3]{PW}.

\begin{lem}
[Comparison Principle] \label{MP}
  Suppose that $u=(u_1,u_2,\cdots,u_k)$ satisfies the following
  uniformly parabolic system of inequalities in $\Om\times(0,T)$.
  \beq
  \left\{
  \bea{c}
  \F{\pa u_1}{\pa t}-\D u_1-\sum\limits_{i=1}^{k}h_{1i}u_i\le0,\\
  \F{\pa u_2}{\pa t}-\D u_2-\sum\limits_{i=1}^{k}h_{2i}u_i\le0,\\
  \vdots\\
  \F{\pa u_k}{\pa t}-\D u_k-\sum\limits_{i=1}^{k}h_{ki}u_i\le0.
  \eea
  \right.
  \eeq
  If $u\le0$ at $t=0$ and on $\pa\Om\times(0,T)$
  and if $h_{ij}$ is bounded and satisfies
  \beq
  h_{ji}\ge0~\mbox{for}~i\neq j, i,j=1,2,\cdots,k,
  \eeq
  then $u\le0$ in $\Om\times(0,T)$.
  Moreover, if there exists $i_0$ such that $u_{i_0}=0$
  at an interior point $(x_0,t_0)$, then $u_{i_0}\equiv0$
  for $t\leq t_0$. Here, we use the notation  $u\le0$ to mean that every component $u_i, i=1,2,\cdots,k$ is nonpositive.
\end{lem}

\begin{prop}
\label{ut vt}
  Let $(u,v)$ satisfies
  \beq
  \label{13}
  \left\{\arraycolsep=1.5pt
  \begin {array}{lll}
  u_{t}-\Delta u=f(x,v)>0,\   \ &\hbox{ in} \,\  \Omega\times(0,T),\\[2mm]
  v_{t}-\Delta v=g(x,u)>0, \   \ &\hbox{ in} \,\  \Omega\times(0,T),\\[2mm]
  u=v=0,\ &\hbox{ on} \ \partial \Omega\times(0,T),\\[2mm]
  u(x,0)=u_0, 
 \ \  v(x,0)=v_0, \ &\hbox{ for}\ x\in\bar{\Omega}.
  \end{array}\right.
  \eeq
Suppose 
${\pa f}/{\pa v}$ and ${\partial g}/{\partial u}$
  being positive and locally bounded. Then if $(u_{0}(x),v_{0}(x))$ is a  subsolution (supersolution) of the corresponding stationary system to \eqref{13},
  there holds $u_t>0,v_t>0$ $(u_t<0,v_t<0)$ in $\Omega$.
\end{prop}

\proof
Consider the case that $(u_0(x),v_0(x))$ is a subsolution of \eqref{13}.
Differentiating system \eqref{13} with respect to $t$ yields
\beq \left\{\arraycolsep=1.5pt
\begin {array}{lll}
(u_t)_{t}-\Delta u_t=\F{\pa f}{\pa v}v_t,\   \ &\hbox{ in} \,\  \Omega\times(0,T),\\[2mm]
(v_t)_{t}-\Delta v_t=\F{\pa g}{\pa u}u_t, \   \ &\hbox{ in} \,\  \Omega\times(0,T),\\[2mm]
u_t=v_t=0,\ &\hbox{ on} \ \partial \Omega\times(0,T),\\[2mm]
u_t(x,0)=\Delta u_0+f(x,v_0)\ge0, \ &\hbox{ for}\ x\in\bar{\Omega},\\[2mm]
v_t(x,0)=\Delta v_0+g(x,u_0)\ge0, \ &\hbox{ for}\ x\in\bar{\Omega}.
\end{array}\right.
\eeq By the maximum principle for parabolic system stated in Lemma
\ref{MP},
we get  $u_t>0$ and $v_t>0$ in $\Omega$.
The  case for $(u_0(x),v_0(x))$ being a supsolution can be proved  similarly.
\proofend

\section{Global existence vs quenching} \label{section Global existence vs quenching}

This section is devoted to prove Theorem \ref{Global existence vs quenching}. Some ideas are borrowed from \cite{fujita}.

\subsection{Proof of Theorem \ref{Global existence vs quenching} $(a)$}
\subsubsection{Proof of Theorem \ref{Global existence vs quenching} $(a_1)$}
For the special case $(u_0,v_0)=(0,0),$
the global existence  can be deduced
directly by sub-super solution method in \cite[Chapter 8]{Pao1}.
Moreover, first the unique global solution
 is bounded by the unique minimal solution
$(w_{\lam,\mu},z_{\lam,\mu})$ of \eqref{E}. Secondly, noting that
$(0,0)$ is a subsolution to \eqref{P},
by Proposition \ref{ut
vt} and assumption \eqref{hhhhh}, we   conclude that $u_t>0$,
$v_t>0$. These imply that
 $(u,v)$ converges as $t\rightarrow+\infty$ to some
functions $\tilde{u}(x),\tilde{v}(x)$ satisfying
 $\tilde{u}\leq w_{\lam,\mu}<1$, $\tilde{v}\leq z_{\lam,\mu}<1$ in
$\Om$.

Let $\varphi(x)\in C^{2}(\bar{\Om})$ and $\varphi|_{\partial\Om}=0$.
Multiplying \eqref{P} by $\varphi$ and integrating over $\Om$, we
arrive at \beq \left\{
\begin{array}{l}
\displaystyle\frac{\mathrm{d}}{\dt}
\displaystyle\int_{\Om}u\varphi\dx-\displaystyle\int_{\Om}u\Delta\varphi\dx
=\displaystyle\int_{\Om}\lam\al(x)\varphi f(v)\dx,\\
\displaystyle\frac{\mathrm{d}}{\dt}
\displaystyle\int_{\Om}v\varphi\dx-\displaystyle\int_{\Om}v\Delta\varphi\dx
=\displaystyle\int_{\Om}\mu \be(x)\varphi g(u)\dx.
\end{array}
\right. \eeq Operating on both sides with
$\displaystyle\frac{1}{T}\displaystyle\int_{0}^{T}$, it follows that
\beq \left\{
\begin{array}{l}
\displaystyle\int_{\Om}\displaystyle\frac{u(x,T)-u_{0}(x)}{T}\varphi\dx
+\displaystyle\int_{\Om}(-\Delta\varphi)\displaystyle\frac{1}{T}\displaystyle\int_{0}^{T}u(x,t)\dt\dx
=\displaystyle\int_{\Om}\lam \al(x)\varphi\displaystyle\frac{1}{T}\displaystyle\int_{0}^{T}f(v)\dt\dx,\\
\displaystyle\int_{\Om}\displaystyle\frac{v(x,T)-v_{0}(x)}{T}\varphi\dx
+\displaystyle\int_{\Om}(-\Delta\varphi)\displaystyle\frac{1}{T}\displaystyle\int_{0}^{T}v(x,t)\dt\dx
=\displaystyle\int_{\Om}\mu
\be(x)\varphi\displaystyle\frac{1}{T}\displaystyle\int_{0}^{T}g(u)\dt\dx.
\end{array}
\right. \eeq Note that \beq
\begin{array}{lll}
\lim\limits_{T\rightarrow+\infty}\displaystyle\frac{u(x,T)-u_{0}(x)}{T}=0,
&\lim\limits_{T\rightarrow+\infty}\displaystyle\frac{v(x,T)-v_{0}(x)}{T}=0,\\
\lim\limits_{T\rightarrow+\infty}\displaystyle\frac{1}{T}\displaystyle\int_{0}^{T}u(x,t)\dx=\tilde{u}(x),
&\lim\limits_{T\rightarrow+\infty}\displaystyle\frac{1}{T}\displaystyle\int_{0}^{T}v(x,t)\dx=\tilde{v}(x),\\
\lim\limits_{T\rightarrow+\infty}\displaystyle\frac{1}{T}\displaystyle\int_{0}^{T}g(u)\dx
=g(\tilde{u}),
&\lim\limits_{T\rightarrow+\infty}\displaystyle\frac{1}{T}\displaystyle\int_{0}^{T}f(v)\dx
=f(\tilde{v}).
\end{array}
\eeq Therefore, by the Lebesgue dominated convergence theorem we get
that as $T\rightarrow+\infty$ \beq \left\{\arraycolsep=1.5pt
\begin {array}{l}
\displaystyle\int_{\Om}\tilde{u}(-\Delta\varphi)\dx=\displaystyle\int_{\Om}\lam \alpha(x)\varphi f(\tilde{v})\dx,\\
\displaystyle\int_{\Om}\tilde{v}(-\Delta\varphi)\dx=\displaystyle\int_{\Om}\mu
\beta(x)\varphi g(\tilde{u})\dx,
\end{array}
\right. \eeq which implies $(\tilde{u},\tilde{v})$ is a weak
solution of \eqref{E}. By the $L^{p}$ estimates of Agmon, Douglis,
Nirenberg \cite{Agmon}, the Sobolev embedding, and the classical
Schauder estimate, we obtain that $(\tilde{u},\tilde{v})$ is a
classical solution of \eqref{E}, and hence
$(\tilde{u},\tilde{v})=(w_{\lam,\mu},z_{\lam,\mu})$.
Since  $u_t>0,v_t>0$ and $w_{\lam,\mu},z_{\lam,\mu}$ are
continuous, by \cite[Theorem 7.13]{Ru} the convergence of the unique
global solution $(u(x,t), v(x,t))$  to
$(w_{\lam,\mu}(x),z_{\lam,\mu}(x))$ is further uniform in $x$, i.e.,
\beq \label{21} \lim_{t\rightarrow
\infty}\big(\|u(x,t)-w_{\lam,\mu}(x)\|_{\infty}+\|v(x,t)-z_{\lam,\mu}(x)\|_{\infty}\big)=0.
\eeq

By comparison principle, the convergence also holds for general case $(u_{0},v_{0})\le(w_{\lam,\mu}, z_{\lam,\mu}).$ Therefore combined with Corollary \ref{converge in C1 norm}, we
complete the proof for Theorem \ref{Global existence vs quenching} $(a_1)$ and Remark \ref{C1 converge}. Furthermore, if $(u_{0}(x),v_{0}(x))$ is further a subsolution of \eqref{E}, by Proposition \ref{ut vt} the convergence is monotone increasing.
 \proofend

\subsubsection{Proof of Theorem \ref{Global existence vs quenching} $(a_2)$}
We first prove Theorem \ref{Global existence vs quenching} $(a_{21})$.
Without loss of generality, we assume that there exists $s\in(0,1)$ such that
 \beq\label{210}
(u_0,v_0)\le \big(sw_1+(1-s)w_{\lam,\mu},sz_1+(1-s)z_{\lam,\mu}\big).
\eeq
Otherwise by Hopf Lemma there exists $t_0>0$ such that $(u(x,t_0),v(x,t_0))$ satisfies the assumption.

For case $(u_0,v_0)= \big(sw_1+(1-s)w_{\lam,\mu},sz_1+(1-s)z_{\lam,\mu}\big)$,
by direct calculations, we can get that $\big(sw_1+(1-s)w_{\lam,\mu},sz_1+(1-s)z_{\lam,\mu}\big)$ is a supersolution of \eqref{E}.
By Proposition \ref{ut
vt} and Proposition \ref{triple}, we  conclude that $u_t<0$,
$v_t<0$ in $\Omega$. Similar to the proof in Theorem \ref{Global existence vs quenching} $(a_1)$, this implies that
 $(u,v)$ converges to $(w_{\lam,\mu},z_{\lam,\mu})$ in $C^1$  as $t\rightarrow+\infty$ and the convergence is monotone decreasing.
 At last, the convergence  holding for general  initial data satisfying $ (0,0)\le(u_0,v_0)\le\not\equiv(w_1,z_1)$ follows directly by the comparison principle. The proof of  Theorem \ref{Global existence vs quenching} $(a_{21})$ and Remark \ref{C1 converge} are completed.

Next, we verify Theorem \ref{Global existence vs quenching} $(a_{22})$. Suppose by contradiction that for some initial value $(u_0,v_0)\ge\not\equiv(w_1,z_1)$,
the solution $(u,v)$ to \eqref{P}  is bounded by $1-\delta$ for all $t\geq0$ and some $0<\delta<1$. Similar to the proof of  Theorem \ref{Global existence vs quenching} $(a_{21})$, without loss of generality, we assume that there exists $\e>0$  such that
$$(\underline{u}_0,\underline{v}_0):=((1+\e)w_1-\e w_{\lam,\mu},(1+\e)z_1-\e z_{\lam,\mu})\leq (u_0,v_0).$$ Secondly, by the convexity of $f$ and $g$ we have
\beq
\D \underline{u}_0+\lam\al(x)f(\underline{v}_0)\ge\not\equiv0,~\D \underline{v}_0+\mu\be(x)g(\underline{u}_0)\ge\not\equiv0,
\eeq
which imply by Proposition \ref{ut vt} that $\underline{u}_t>0,\underline{v}_t>0$ in $\Omega$. Here $(\underline{u},\underline{v})$ is the solution to \eqref{P} with $(\underline{u}(x,0),\underline{v}(x,0))=(\underline{u}_0,\underline{v}_0)$. Thus $(\underline{u},\underline{v})$   converges increasingly to some
$(\bar{w},\bar{z})$ and hence $(u,v)$ converges to $(\bar{w},\bar{z})$ by the comparison principle. Note that it can be checked similarly to Theorem \ref{Global existence vs quenching} $(a_1)$ that $(\bar{w},\bar{z})$ is a classical solution to \eqref{E}. We  have $(w_{\lam,\mu},z_{\lam,\mu})\ll (w_1,z_1)\ll(\bar{w},\bar{z}),$ which contradicts Proposition \ref{triple}. The proof of  Theorem \ref{Global existence vs quenching} $(a_{22})$ is completed.

\subsection{Proof of Theorem \ref{Global existence vs quenching} $(b)$}
\label{Touchdown behavior}

 We will only  prove the case $(u_{0},v_{0})\equiv(0,0)$, then
the holding for general nonnegative initial data follows directly by
comparison principle.

Let $(\lam,\mu)\in\mathcal {O}_2$.  Suppose on the contrary that the
local solution $(u,v)$ (see Theorem B) exists globally, i.e. $0\le
u<1$, $0\le v<1$ for all $t\ge0$. Take $\delta>1$,
$a={\lam}/{\delta}$, $b={\mu}/{\delta}$. Since
$U={u}/{\delta}<u$, $V={v}/{\delta}<v$, it then indicates that
$U\le{1}/{\delta}<1$, $V\le{1}/{\delta}<1$, and by the monotone
increasing of $f,g$ there holds \beq \left\{\arraycolsep=1.5pt
\begin {array}{lll}
U_t-\Delta U={\lam\al(x)f(v)}/{\delta}\ge a\al(x)f(V),\   \ &\hbox{ in} \,\ \Omega\times(0,T),\\[2mm]
V_t-\Delta V={\mu \beta(x)g(u)}/{\delta}\ge b \beta(x)g(U), \   \ &\hbox{ in} \,\  \Omega\times(0,T),\\[2mm]
U=V=0,\ &\hbox{ on} \ \partial \Omega\times(0,T),\\[2mm]
U(x,0)=0,\ V(x,0)=0, \ &\hbox{ for}\ x\in\bar{\Omega}.
\end{array}
\right. \eeq Hence $(U,V)$ is a supersolution of \beq \label{(w,z)}
\left\{\arraycolsep=1.5pt
\begin {array}{lll}
\mathcal{U}_t-\Delta\mathcal{U}=a\alpha(x)f(\mathcal{V}),\   \ &\hbox{ in} \,\  \Omega\times(0,T),\\[2mm]
\mathcal{V}_t-\Delta\mathcal{V}=b \beta(x)g(\mathcal{U}), \   \ &\hbox{ in} \,\  \Omega\times(0,T),\\[2mm]
\mathcal{U}=\mathcal{V}=0,\ &\hbox{ on} \ \partial \Omega\times(0,T),\\[2mm]
\mathcal{U}(x,0)=0,\ \mathcal{V}(x,0)=0, \ &\hbox{ for}\
x\in\bar{\Omega}.
\end{array}
\right. \eeq Therefore \eqref{(w,z)} has a global classical solution
$(\mathcal{U}(x,t),\mathcal{V}(x,t))$, since $0\le\mathcal{U}\le
U\le{1}/{\delta}<1$, $0\le\mathcal{V}\le V\le{1}/{\delta}<1$. Note
that there further holds \beq \label{L2 H2 norm} \left\{
\begin{array}{l}
\lim\limits_{t\rightarrow+\infty}(\|\mathcal{U}_t\|_2^2+\|\mathcal{V}_t\|_2^2)=0,\\
\sup\limits_{t>1}(\|\mathcal{U}\|_{H^2(\Om)}+\|\mathcal{V}\|_{H^2(\Om)})<+\infty,
\end{array}
\right. \eeq which can be proved similarly to Proposition \ref{22}.
By Sobolev embedding theorem, one can have that there exists a
subsequence $\{t_{j}\}_{j=1}^{\infty}$ such that
$t_{j}\rightarrow+\infty$,
$(\mathcal{U}(\cdot,t_{j}),\mathcal{V}(\cdot,t_{j})$ converges
strongly in $H^1_0(\Om)$ to
$(\mathcal{U}_{\infty},\mathcal{V}_{\infty})$. Now take $\phi\in
H_0^1(\Om)$. Multiplying \eqref{(w,z)} by $\phi$ and integrating by
parts with respect to $x$ yields,  \beq \left\{\arraycolsep=1.5pt
\begin {array}{l}
\dis\int_{\Om}\phi\mathcal{U}_t(\cdot,t_{j})\dx
+\dis\int_{\Om}\nabla\mathcal{U}(\cdot,t_{j})\nabla\phi\dx
=\dis\int_{\Om}a\alpha(x)\phi f(\mathcal{V}(\cdot,t_{j}))\dx,\\[2mm]
\dis\int_{\Om}\phi\mathcal{V}_t(\cdot,t_{j})\dx
+\dis\int_{\Om}\nabla\mathcal{V}(\cdot,t_{j})\nabla\phi\dx
=\dis\int_{\Om}b\beta(x)\phi g(\mathcal{U}(\cdot,t_{j}))\dx.
\end{array}
\right. \eeq Let $t_{j}\rightarrow+\infty$ and  we will get
 that $(\mathcal{U}_{\infty},\mathcal{V}_{\infty})$ is a weak
solution of \beq \left\{\arraycolsep=1.5pt
\begin {array}{lll}
-\Delta w=a \alpha(x) f(z),\ \ &\hbox{in} \,\  \Omega,\\[2mm]
-\Delta z=b \beta(x) g(w), \ \ &\hbox{in} \,\  \Omega,\\[2mm]
u=v=0,\ &\hbox{on} \ \partial \Omega.
\end{array}
\right. \eeq Chose  $\delta$  close to 1 such that $(a,b)\in\mathcal
{O}_2$. Then we   get a contradiction with Theorem A. The proof of
this proposition is therefore completed. \proofend

\subsection{Proof of Theorem \ref{Global existence vs quenching} $(c)$}

\label{Quenching for large initial data}

First note from Theorem B that the solution of \eqref{P} exists locally.
Suppose that $(u,v)$ is the solution of \eqref{P} for $0\le t<T$. Then there holds $0\leq u,v<1$ for $0\le t<T$.
Define
$F(t)=\dis\int_{\Om}u\phi\dx$. Then some calculation with
integrating by parts leads to \beq \label{F'(t)}
\begin{array}{l}
F'(t)=\dis\int_{\Om}u_t\phi\dx
=\dis\int_{\Om}(\Delta u+\lam\al(x)f(v))\phi\dx
=-\lam_{1}\dis\int_{\Om}u\phi\dx+\lam\dis\int_{\Om}\al(x)\phi
f(v)\dx\ \text{for}\ t\in[0,T).
\end{array}
\eeq
Note that there holds \beq
F(t)=\dis\int_{\Om}u\phi\dx
\le\dis\int_{\Om}\F{\sqrt{f(v)}\sqrt{\phi}\sqrt{\phi}}{\sqrt{f(0)}}\dx
\le\bigg(\dis\int_{\Om}\F{\al f(v)\phi}{f(0)}\dx\bigg)^{\frac12}
\bigg(\dis\int_{\Om}{\phi}/{\al}\dx\bigg)^{\frac12}
\ \text{for}\ t\in[0,T),\eeq
by \eqref{hhhhh} and $u<1.$
Consequently,
$F^2(t)\le\bigg(\dis\int_{\Om}\F{\al f(v)\phi}{f(0)}\dx\bigg)\bigg(\dis\int_{\Om}{\phi}/{\al}\dx\bigg)$. Combining this
inequality with \eqref{F'(t)} yields
\beq
F'(t)+\lam_{1}F(t)\ge\big({\lam f(0)F^2(t)}\big)/\bigg({\dis\int_{\Om}{\phi}/{\al}\dx}\bigg). \eeq
Thus
$\F{\dis \mathrm{d}}{\dis \dt}\bigg(-\F{1}{e^{\lam_1 t}F(t)}\bigg)
\ge\big({\lam f(0)e^{-\lam_1 t}}\big)/\bigg(
{\dis\int_{\Om}{\phi}/{\al}\dx}\bigg)$.
By integrating on both side there hold
\beq
e^{-\lam_1t}\bigg(-\F{1}{F(t)}\bigg)+\F{1}{F(0)}
\ge\F{\lam f(0)}{\lam_1\dis\int_{\Om}{\phi}/{\al}\dx}(1-e^{-\lam_1t}),
\eeq
and then
\beq \label{11113}
e^{-\lam_1t}\bigg(\F{\lam f(0)}{\lam_1\dis\int_{\Om}{\phi}/{\al}\dx}-\F{1}{F(t)}\bigg)
\ge\F{\lam f(0)}{\lam_1\dis\int_{\Om}{\phi}/{\al}\dx}-\F{1}{F(0)}.
\eeq
Note that $$\F{\lam f(0)}{\lam_1\dis\int_{\Om}{\phi}/{\al}\dx}-\F{1}{F(t)}\ge
e^{\lam_1t}\bigg(\F{\lam f(0)}{\lam_1\dis\int_{\Om}{\phi}/{\al}\dx}-\F{1}{F(0)}\bigg)\ge0$$
by
$F(0)=\dis\int_{\Om}u_{0}\phi\dx>({\lam_{1}\dis\int_{\Om}{\phi}/{\al}\dx})/({\lam f(0)}).$
From \eqref{11113} it can be deduced that
\beq
e^{-\lam_1t}\ge\F{\lam f(0)-\bigg({\lam_1\dis\int_{\Om}{\phi}/{\al}\dx}\bigg)/{F(0)}}{\lam f(0)-\bigg({\lam_1\dis\int_{\Om}{\phi}/{\al}\dx}\bigg)/{F(t)}},
\eeq
which implies that
\beq
t\le\F{1}{\lam_1}\ln{\F
{\lam f(0)-\bigg({\lam_1\dis\int_{\Om}{\phi}/{\al}\dx}\bigg)/{F(t)}}
{\lam f(0)-\bigg({\lam_1\dis\int_{\Om}{\phi}/{\al}\dx}\bigg)/{F(0)}}}
\le\F{1}{\lam_1}\ln{\F
{\lam f(0)-{\lam_1\dis\int_{\Om}{\phi}/{\al}\dx}}
{\lam f(0)-\bigg({\lam_1\dis\int_{\Om}{\phi}/{\al}\dx}\bigg)/{F(0)}}}\ \text{for any}\ t\in[0,T).
\eeq
Hence $u$ must quench at finite time $T^*_u$ and $T^*_u\le\F{1}{\lam_1}\ln{\F
{\lam f(0)-{\lam_1\dis\int_{\Om}{\phi}/{\al}\dx}}
{\lam f(0)-\bigg({\lam_1\dis\int_{\Om}{\phi}/{\al}\dx}\bigg)/{F(0)}}}$.

Similarly, for any given $\mu>0$, once $\dis\int_{\Om}v_{0}\phi\dx>{\lam_{1}\dis\int_{\Om}{\phi}/{\be}\dx}/
{\mu g(0)}$, $v(x,t)$
  must quench at a finite time $T_v^*$ and $T^*_v\le\F{1}{\lam_1}\ln{\F{\mu g(0)-\lam_1\dis\int_{\Om}{\phi}/{\be}\dx}{\mu g(0)-\bigg(\lam_1\dis\int_{\Om}{\phi}/{\be}\dx\bigg)/G(0)}}$. The proof is therefore completed.
\proofend

\section{Convergence rate}\label{convergence rate}
\label{Convergence rate for MEMS system}

This section is devoted to prove Theorem \ref{Convergence rate}. Here, we point out that compared to single MEMS equation, some new ideas and
techniques are needed to obtain the convergence rate.

 To obtain
the convergence rate, we need to consider the stability
of
$(w_{\lam,\mu},z_{\lam,\mu})$. For this purpose, we first 
show a related lemma as follows.

\begin{lem}
\label{aaa6} The problem\beq \label{7}
 \left\{\arraycolsep=1.5pt
\begin {array}{lll}
-\Delta \varphi-\lam\al(x)f'(z_{\lam,\mu})\psi=\nu\varphi,\ \ &\hbox{in} \,\  \Omega,\\[2mm]
-\Delta \psi-\mu\be(x)g'(w_{\lam,\mu})\varphi=\nu\psi, \ \ &\hbox{in} \,\  \Omega,\\[2mm]
\varphi=\psi=0,\ &\hbox{on} \ \partial \Omega.
\end{array}\right.
\eeq  has a first eigenvalue $\nu_1>0$ (which means the minimal
solution $(w_{\lam,\mu},z_{\lam,\mu})$ is stable) with strictly
positive eigenfunction $(\varphi_1,\psi_1)$, that is, $\varphi_1>0,
\psi_1>0$ in $\Omega.$ Moreover $\varphi_1$ and $\psi_1$ are smooth.
\end{lem}
This  result is standard. For the proof, see e.g., \cite[Theorem
1.5]{MON} and \cite[p10]{Clemente}.

Next, before verifying Theorem \ref{Convergence rate}, we shall introduce the following two  lemmas and a
proposition.

\begin{lem}
\label{16}
Given a smooth bounded domain $\Omega$ in $\mathbb{R}^N$.  Suppose 
  $a(x,t)\in C([0,+\infty),C^1(\bar{\Omega}))$, $a\geq0$ in $\bar{\Omega}\times [0,+\infty),$ $b(x)\in C^1(\bar{\Omega}),$ $b>0$ in $\Omega,$
  $a=b=0$ on $\partial\Om$,
  $\F{\pa b}{\pa\vec{n}}<0$ on $\partial\Om$,
   $\lim\limits_{t\rightarrow+\infty}\|a(\cdot,t)\|_{C^1}=0$.
  Then there exists $T_0>0$ such that
  $a(x,t)\le b(x)$ in $\Om$ for all  $t>T_0$.
  Here $\vec{n}$ denotes the outward unit normal vector on $\partial \Omega$.
\end{lem}
\proof Since $\F{\pa b}{\pa\vec{n}}|_{\pa\Om}<0$ and $b(x)\in
C^1(\bar{\Omega}),$ there exists a constant $\varepsilon>0$ such
that for all $x\in\Omega_{\varepsilon}:=\{x\in
\bar{\Omega}|dist(x,\partial \Omega)\leq \varepsilon\}$, there holds
$b(x)=b(x)-b(x_0)\ge C_0|x-x_0|,$ where $x_0\in \partial \Omega$
satisfying $(x-x_0)\parallel\vec{n}$ and $C_0>0$ is a constant
independent on $x.$ On the other hand, for all
$x\in\Omega_{\varepsilon}$, there also holds
$a(x,t)=a(x,t)-a(x_0,t)\le    \|a(\cdot,t)\|_{C^1} |x-x_0|.$ Note
that
 $\lim\limits_{t\rightarrow+\infty}\|a(\cdot,t)\|_{C^1}=0.$ Therefore, there holds $ \|a(\cdot,t)\|_{C^1} \leq C_0$ for $t$ large enough and
 it follows that $a(x,t)\leq b(x)$ on $\Omega_{\varepsilon}$ for $t$ large enough. At last,
 it is obviously that for any given subset $\mathring{\Omega}\subset\Omega,$   $a(x,t)\le b(x)$ on $\mathring{\Omega}$ for $t$ large enough. Hence, we conclude this lemma.
 \proofend

\begin{lem}
\label{17}
  For the solution $(u,v)$ to Problem \eqref{P} with $(\lambda,\mu)\in\mathcal {O}_1$,
  if $(u_{0}(x),v_{0}(x))$ is a  subsolution (supersolution) of \eqref{E},
  then there exist $c_1,c_2\in \mathbb{R}^+$ such that
\beq \label{a5} u_t\geq c_1 v_t\geq 0\  \text{and}\  v_t\geq c_2 u_t\geq 0\quad (u_t\leq c_1 v_t\leq 0\  \text{and}\  v_t\leq c_2 u_t\leq 0).
\eeq
\end{lem}

\proof
We only show the proof for $(u_{0}(x),v_{0}(x))$ being a  subsolution, since the proof for $(u_{0}(x),v_{0}(x))$ being a supersolution is totally similar.
  Let
  \beq
  U=u_t-c_1 v_t,\quad V=v_t-c_2u_t.
  \eeq
  Note by Proposition \ref{ut vt} that $u_t\geq0$, $v_t\geq0$.
  It can be deduced that
  \beq
  \left\{
  \bea{l}
  U_t-\Delta U+c_1\mu \beta(x)g'(u)U=(\lam\al f'(v)-c_1^2\mu\be g'(u))v_t,\\
  U|_{\partial \Omega}=0,\\
  U(x,0)=\Delta u_0+\lambda \alpha(x)f(v_0)-c_1(\Delta v_0+\mu\beta(x)g(u_0)).
  \eea
  \right.
  \eeq
    Applying comparison principle, we have that $u_t-c_1 v_t=U\geq 0$ provided that
    \beq
    c_1\leq \min\bigg\{\sqrt{\frac{\lambda}{\mu}\inf_{\Omega}\frac{\alpha(x)}{\beta(x)}\frac{f'(0)}{g'(\|w_{\lam,\mu}\|_{\infty})}},
    \inf_{\Omega}\frac{\Delta u_0+\lambda \alpha(x)f(v_0)}{\Delta v_0+\mu \beta(x)g(u_0)}\bigg\}.
    \eeq
    Here $0 \leq u,v<1$, the nonnegativity of
    $\lam,\mu, \al, \be, u_t, v_t, f'(s), g'(s)\ \text {for}\ 0\leq s<1$
    and the monotonicity of $f'(s), g'(s)$ are used.
    Similarly, it can be proved that $v_t-c_2 u_t\geq 0$ provided
   \beq
    c_2\leq \min\bigg\{\sqrt{\frac{\mu}{\lambda}\inf_{\Omega}\frac{\beta(x)}{\alpha(x)}\frac{g'(0)}{f'(\|z_{\lam,\mu}\|_{\infty})}},
    \ \inf_{\Omega}\frac{\Delta v_0+\mu \beta(x)g(u_0)}{\Delta u_0+\lambda \alpha(x)f(v_0)}\bigg\}.
    \eeq
    This completes the proof of \eqref{a5}.
\proofend

Without causing confusion, for simplicity  we use $(w,z)$ instead of
$(w_{\lam,\mu},z_{\lam,\mu})$ to denote the  minimal solution of
problem \eqref{E} in the rest part of this subsection.
\begin{prop}\label{22}
Assume that $(u_0,v_0)\leq(w,z)$ and is a   subsolution of \eqref{E}, or $(w,z)\leq(u_0,v_0)\le\not\equiv
(w_1, z_1)$(if $(w_1,z_1)$ exists) and is a   supsolution of \eqref{E}.
  Let $(u,v)$ be the unique global solution of \eqref{P},
  then we have
\beq \label{11}
  \lim_{t\rightarrow +\infty}\|u_t\|_2= \lim_{t\rightarrow +\infty}\|v_t\|_2=0
\eeq and \beq\label{12}
  \|u\|_{H^3}+  \|v\|_{H^3}\leq C(\delta),~
  \mbox{for all}~t\ge\delta>0.
\eeq
\end{prop}

\proof First we claim that \beq \label{a2}
  \|\nabla u\|_2+\|\nabla v\|_2\leq  C(u_0,v_0,w,z),
  \eeq
  where $C$ is a constant independent of time $t$.
To prove this claim, we denote $\xi=u-w, \eta=v-z$. Then it follows
from system \eqref{P} and \eqref{E} that
 \beq\label{xi} \left\{\arraycolsep=1.5pt
\begin {array}{lll}
\xi_t-\Delta \xi=\lam\al(x)(f(v)-f(z)),\   \ &\hbox{ in} \,\  \Omega\times(0,T),\\[2mm]
\eta_t-\Delta \eta=\mu\be(x)(g(u)-g(w)), \   \ &\hbox{ in} \,\  \Omega\times(0,T),\\[2mm]
\xi=\eta=0,\ &\hbox{ on} \ \partial \Omega\times(0,T),\\[2mm]
\xi(x,0)=u_0(x)-w(x),\ \eta(x,0)=v_0(x)-z(x), \ &\hbox{ for}\
x\in\bar{\Omega}.
\end{array}\right.
\eeq Multiplying the first equation of  \eqref{xi} by $\xi_t$ yields
that \beq \frac{1}{2}\frac{d}{\dt}\|\nabla \xi\|_2^2+\|\xi_t\|_2^2
=\lam\int_{\Omega}\al(x)[f(v)-f(z)]\xi_t\dx\leq0, \eeq
where assumption \eqref{hhhhh} and Proposition \ref{ut vt} are used.
The above inequality then implies that \beq \frac{d}{dt}\|\nabla
\xi\|_2^2\leq0 \eeq and hence \beq \|\nabla \xi\|_2\leq\|\nabla
\xi(x,0)\|_2\leq C(u_0,v_0,w,z). \eeq $\|\nabla \eta\|_2\leq
C(u_0,v_0,w,z)$ can be obtained similarly. Then \eqref{a2} follows
by \beq
  \|\nabla u\|_2+\|\nabla v\|_2\leq \|\nabla\xi\|_2+\|\nabla w\|_2+\|\nabla \eta\|_2+\|\nabla z\|_2\leq  C(u_0,v_0,w,z).
  \eeq

Next, we show that \beq\label{a6}
\int_0^{+\infty}(\|u_t\|_2^2+\|v_t\|_2^2)\dt\leq C. \eeq After
multiplying  equations in \eqref{P} by $v_t$ and $u_t,$
respectively, adding them up and integrating over $\Omega$, we can
see that Problem \eqref{P} admits a Lyapunov function \beq
E(u,v)=\int_{\Omega}\big(\nabla u\nabla
v-\mathbb{F}(x,v)-\mathbb{G}(x,u)\big)\dx, \eeq where
$\mathbb{F}(x,v)=\lambda \alpha(x)\dis\int_{0}^{v}f(s)\ds$,
$\mathbb{G}(x,u)=\mu \beta(x)\dis\int_{0}^{u}g(s)\ds,$  and there
holds
\beq\label{a4} \frac{\mathrm{d}}{\dt}E(u,v)+2\int_{\Omega}u_tv_t\dx=0.
\eeq Note that $0<w<1, 0<z<1$. By assumption \eqref{hhhhh} and
\eqref{a2}, integrating \eqref{a4} with respect to $t$ yields


\beq\label{aaa3}
\bea{l}
2\dis\int_0^{+\infty}\int_{\Omega}u_tv_t\dx\mathrm{d}\tau \leq
E(u_0,v_0)+\bigg|\int_{\Omega}\nabla u\nabla v\dx\bigg|
+\int_{\Omega}\mathbb{F}(x,C_1)\dx
+\dis\int_{\Omega}\mathbb{G}(x,C_2)dx\leq
C, \eea\eeq
where $C_1=\max_{\Omega}z$ and $C_2=\max_{\Omega}w$ for  $(u_0,v_0)\leq(w,z)$, $C_1=\max_{\Omega}z_1$ and $C_2=\max_{\Omega}w_1$ for $(w,z)\leq(u_0,v_0)\le\not\equiv
(w_1, z_1)$(if $(w_1,z_1)$ exists).
Then
 \eqref{a6} can be concluded by \eqref{aaa3}, Proposition \ref{ut vt} and Lemma \ref{17}.

Differentiating the first equation in \eqref{P} with respect to $t$
yields \beq\label{a9} u_{tt}-\Delta u_t=\lam\al(x)f'(v)v_t. \eeq
Multiplying \eqref{a9} by $u_t$ and integrating over $\Omega$, by
Lemma \ref{17} and assumption \eqref{hhhhh} we have \beq\label{a8}
\frac{1}{2}\frac{\mathrm{d}}{\dt}\|u_t\|_2^2+\|\nabla
u_t\|_2^2=\int_{\Omega}\lam\al f'(v)v_t u_t \dx\leq C\|u_t\|_2^2.
\eeq By  Young's inequality we get \beq \frac{\mathrm{d}}{\dt}\|u_t\|_2^2\leq
C_1\|u_t\|_2^4+C_2. \eeq
 Then by \eqref{a6} and \cite[Lemma 6.2.1]{Zheng},
 we get $\lim\limits_{t\rightarrow +\infty}\|u_t\|_2=0$,
 while $\lim\limits_{t\rightarrow +\infty}\|v_t\|_2=0$ can be obtained similarly and \eqref{11} follows.

Now integrating \eqref{a8} with respect to $t,$ by \eqref{a6} we
obtain \beq \frac{1}{2}\|u_t\|_2^2+\int_0^{t}\|\nabla
u_t\|_2^2\mathrm{d}\tau\leq
\frac{1}{2}\|u_t(0)\|_2^2+C\int_0^{+\infty}\|u_t\|_2^2\mathrm{d}\tau\leq C,
\eeq which implies obviously \beq \label{14} \int_0^{t}\|\nabla
u_t\|_2^2\mathrm{d}\tau\leq C. \eeq Multiplying \eqref{a9} by $-\Delta u_t$
and integrating over $\Omega$, by Lemma \ref{17} we have \beq
\frac{1}{2}\frac{\mathrm{d}}{\dt}\|\nabla u_t\|_2^2+\|\Delta u_t\|_2^2
=\lam\int_{\Omega}\al f'(v)v_t(-\Delta u_t)\dx\leq C\|u_t\|_2\|\Delta
u_t\|_2\leq C\|u_t\|_2^2+\frac{1}{2}\|\Delta u_t\|_2^2, \eeq which
yields \beq\label{a10} \frac{\mathrm{d}}{\dt}\|\nabla u_t\|_2^2+\|\Delta
u_t\|_2^2\leq C\|u_t\|_2^2. \eeq Multiplying \eqref{a10} by $t$,
then integrating with respect to $t$ in $[0,t]$, by \eqref{a6} and
\eqref{14} there holds \beq t\|\nabla u_t\|_2^2+\int_0^t\tau
\|\Delta u_t\|_2^2\mathrm{d}\tau\leq \int_0^t\|\nabla
u_t\|_2^2\mathrm{d}\tau+Ct\int_0^t\|u_t\|_2^2\mathrm{d}\tau\leq C_1+C_2t. \eeq Thus,
for $t\geq\delta>0,$ we have \beq \label{aaa4} \|\nabla
u_t\|_2^2\leq\frac{C_1}{t}+C_2\leq \frac{C_1}{\delta}+C_2 \eeq and
it follows \beq\label{aaa5} \|u_t\|_{H^1}\leq
C(\delta)~\mbox{for}~t\ge\delta. \eeq Now we can deduce from the
equation in \eqref{P} and the regularity theory for the elliptic
problem (see e.g. \cite{Gilbarg,Zheng}) \beq
 \left\{\arraycolsep=1.5pt
\begin {array}{lll}
-\Delta u=\lam\al(x)f(v)-u_t,\ \ &\hbox{in} \,\  \Omega,\\[2mm]
u=0,\ &\hbox{on} \ \partial \Omega,
\end{array}\right.
\eeq that \beq \label{15} \|u(\cdot,t)\|_{H^3}\leq
C(\|f(v)\|_{H^1}+\|u_t\|_{H^1})\leq C(\delta)(1+\|f'(v)\nabla
v\|_{2})\leq C(\delta)(1+C\|\nabla v\|_2)\leq C(\delta) \eeq by
\eqref{aaa5} and \eqref{a2}. $\|v(\cdot,t)\|_{H^3}$ can be treated
similarly. In conclusion, we obtain \eqref{12}. The proof of Proposition \ref{22} is completed.
\proofend

\begin{cor}
\label{converge in C1 norm}  Suppose that the conditions in Proposition \ref{22}  are satisfied.
 For $1\leq n\leq 3,$ there holds
  \beq
  \label{a11}
  \lim_{t\rightarrow\infty}\|\xi(\cdot,t)\|_{C^1}
  =\lim_{t\rightarrow\infty}\|u(\cdot,t)-w\|_{C^1}=0,\ \lim_{t\rightarrow\infty}\|\eta(t)\|_{C^1}=\lim_{t\rightarrow\infty}\|v(t)-z\|_{C^1}=0.
  \eeq
\end{cor}
\proof
 Note by \eqref{12} that $\xi,\eta\in H^3(\Om)$,
and $H^3(\Om)\hookrightarrow\hookrightarrow C^1(\Om)$ for $1\le
n\leq 3$ by Sobolev compact embedding theorem. Thanks to \eqref{21},
\eqref{a11} follows by the relative compactness of $\xi(t),\eta(t)$
in $C^1$ and the uniqueness of the limits.
\proofend

Now we present the proof of Theorem \ref{Convergence rate} as follows.\\

\proofconrate
Multiplying  equations in \eqref{xi} by $\xi$ and $\eta,$
respectively, adding them up and integrating over $\Omega$ yields
\beq\label{9} \frac{\mathrm{d}}{\dt}\int_{\Omega}\bigg(\frac{1}{2}\xi^2
+\frac{1}{2}\eta^2\bigg)\dx+\|\nabla \xi\|_2^2+\|\nabla \eta\|_2^2
=\int_{\Omega}\bigg(\lam\al[f(z)-f(v)](-\xi)+\mu\be[g(w)-g(u)](-\eta)\bigg)\dx.
\eeq
Rewrite  equations in \eqref{xi} as \beq\label{xi1}
\left\{\arraycolsep=1.5pt
\begin {array}{lll}
\xi_t-\Delta \xi-\lam\al f'(z)\eta=\lam\al(f(v)-f(z)-f'(z)\eta),\   \ &\hbox{ in} \,\  \Omega\times(0,T),\\[2mm]
\eta_t-\Delta \eta-\mu\be g'(w)\xi=\mu\be(g(u)-g(w)-g'(w)\xi), \   \
&\hbox{ in} \,\  \Omega\times(0,T).
\end{array}\right.
\eeq 
By the
convexity of $f$ and $g$, it is easy to deduce that
$f(v)-f(z)-f'(z)\eta\geq 0$ and $g(u)-g(w)-g'(w)\xi\geq 0.$ Thus it
follows that \beq\label{xi2} \left\{\arraycolsep=1.5pt
\begin {array}{lll}
\xi_t-\Delta \xi-\lam\al f'(z)\eta\geq 0,\   \ &\hbox{ in} \,\  \Omega\times(0,T),\\[2mm]
\eta_t-\Delta \eta-\mu\be g'(w)\xi\geq 0, \   \ &\hbox{ in} \,\
\Omega\times(0,T).
\end{array}\right.
\eeq\\

\textbf{Step 1: Consider $(u_0,v_0)=(0,0)$.}

Let's keep in mind that that $(0,0)$ is a subsolution to \eqref{E}. Multiplying  inequalities in \eqref{xi2} by $\psi_1$ and
$\varphi_1,$ respectively, adding them up and integrating over
$\Omega$ yields
\beq\label{xi3}
\dis\int_{\Omega}(\psi_1\xi+\varphi_1\eta)_t\dx+\nu_1\dis\int_{\Omega}(\psi_1\xi+\varphi_1\eta)\dx
\geq 0. \eeq Here, $\nu_1$ is the principal eigenvalue of problem
\eqref{7} and  $(\varphi_1,\psi_1)$ is  the  corresponding positive
eigenfunction.  Multiplying \eqref{xi3} by $-1$, then adding it to
\eqref{9} yields \beq\label{a13}
\begin{split}
&\frac{\mathrm{d}}{\dt}\int_{\Omega}\bigg(\frac{1}{2}\xi^2+\frac{1}{2}\eta^2+\big[\psi_1(-\xi)+\varphi_1(-\eta)\big]\bigg)\dx+\|\nabla\xi\|_2^2
+\|\nabla \eta\|_2^2+\nu_1\int_{\Omega}\big[\psi_1(-\xi)+\varphi_1(-\eta)\big]\dx\\
\leq&
\int_{\Omega}\bigg(\lam\al[f(z)-f(v)](-\xi)+\mu\be[g(w)-g(u)](-\eta)\bigg)\dx.
\end{split}
\eeq Now we claim that there exists $T_0>0$ such that for any
$t>T_0,$ there holds \beq\label{a12} f(z)-f(v)\leq
\frac{\nu_1}{2\lambda\|\alpha\|_{\infty}}\psi_1, \quad g(w)-g(u)\leq
\frac{\nu_1}{2\mu\|\beta\|_{\infty}}\varphi_1. \eeq In fact,
recalling \eqref{7}, by Lemma \ref{aaa6} we have \beq \label{8}
 \left\{\arraycolsep=1.5pt
\begin {array}{lll}
-\Delta \varphi_1=f'(z)\psi_1+\nu_1\varphi_1\geq 0,\ \ &\hbox{in} \,\  \Omega,\\[2mm]
-\Delta \psi_1=g'(w)\varphi_1+\nu_1\psi_1\geq0, \ \ &\hbox{in} \,\  \Omega,\\[2mm]
\varphi_1=\psi_1=0,\ &\hbox{on} \ \partial \Omega.
\end{array}\right.
\eeq Thus, by Hopf lemma there holds $-\frac{\partial
\varphi_1}{\partial \vec{n}}\geq \varepsilon_0,$ $-\frac{\partial
\psi_1}{\partial \vec{n}}\geq \varepsilon_0$ on $\partial \Omega$
for some $\varepsilon_0>0.$
Then \eqref{a12} follows by Lemma \ref{aaa6}, Lemma \ref{16} and
Corollary \ref{converge in C1 norm}.

Combining \eqref{a13},\eqref{a12} with the Poincar$\acute{e}$
inequality $\|u\|_2\leq \frac{1}{\lambda_1}\|\nabla u\|_2$ for any
$u\in H_0^1(\Omega)$ with $\lambda_1>0$ being the first eigenvalue
of $-\Delta$ on $H_0^1(\Omega),$ we get \beq\label{a14}
\begin{split}
&\frac{\mathrm{d}}{\dt}\int_{\Omega}\bigg(\frac{1}{2}\xi^2+\frac{1}{2}\eta^2+\big[\psi_1(-\xi)+\varphi_1(-\eta)\big]\bigg)dx+\lambda_1\|\xi\|_2^2
+\lambda_1\|\eta\|_2^2+\frac{\nu_1}{2}\int_{\Omega}\big[\psi_1(-\xi)+\varphi_1(-\eta)\big]\dx
\leq0.
\end{split}
\eeq Let
$Y=\dis\int_{\Omega}(\xi^2+\eta^2+2\big[\psi_1(-\xi)+\varphi_1(-\eta)\big])\dx.$
Note that \beq\label{aaa7}
\int_{\Omega}\big[\psi_1(-\xi)+\varphi_1(-\eta)\big]\dx\geq 0. \eeq
By \eqref{a14} there holds \beq \frac{\mathrm{d}Y}{dt}+\gamma Y\leq 0, \quad
\gamma=\min\bigg\{2\lambda_1,\frac{\nu_1}{2}\bigg\}, \eeq which
yields $Y\leq Y(0)e^{-\gamma t}.$ Then by noting \eqref{aaa7} again
it follows
 that \beq \label{215}
  \|u(x,t)-w(x)\|_2^2+\|v(x,t)-z(x)\|_2^2\leq Y(t)
  \leq C_0 \exp\bigg({-\min\bigg\{\lambda_1,\frac{\nu_1}{2}\bigg\}t}\bigg),\quad \text{for}\ t>T_0
\eeq
with $C_0=Y(0)=\|w_{\lam,\mu}\|_2^2+\|z_{\lam,\mu}\|_2^2
+2\|\psi_1w_{\lam,\mu}+\varphi_1z_{\lam,\mu}\|_1.$\\

\textbf{Step 2: Consider $(u_0,v_0)= \big(sw_1+(1-s)w_{\lam,\mu},sz_1+(1-s)z_{\lam,\mu}\big)$ for some $0<s<1$.}

Let's keep in mind that that this $(u_0,v_0)$ is a supsolution to \eqref{E}.
Combining  \eqref{9} with \eqref{xi1} yields
\beq\label{226}
\begin{array}{ll}
&\F{\mathrm{d}}{\dt}\int_{\Omega}\bigg(\frac{1}{2}\xi^2+\frac{1}{2}\eta^2+\big[\psi_1\xi+\varphi_1\eta\big]\bigg)\dx+\|\nabla\xi\|_2^2
+\|\nabla \eta\|_2^2+\nu_1\int_{\Omega}\big[\psi_1\xi+\varphi_1\eta\big]\dx\\
=&
\dis\int_{\Omega}\bigg(\lam\al[f(v)-f(z)]\xi+\mu\be[g(u)-g(w)]\eta\bigg)\dx
+\dis\int_{\Omega}\lam\al(f(v)-f(z)-f'(z)\eta)\psi_1\dx\\
&+\dis\int_{\Omega}\mu\be(g(u)-g(w)-f'(w)\xi)\varphi_1\dx.
\end{array}
\eeq
Similar to \eqref{a12}, there hold
\beq \label{227}
f(v)-f(z)\leq
\frac{\nu_1}{4\lambda\|\alpha\|_{\infty}}\psi_1, \quad g(u)-g(w)\leq
\frac{\nu_1}{4\mu\|\beta\|_{\infty}}\varphi_1\eeq
for $t$ sufficiently large.
Meanwhile, by Corollary \ref{converge in C1 norm} one can obtain that for $t$ sufficiently large,
\beq\label{228}
\lam\al(f(v)-f(z)-f'(z)\eta)\psi_1\leq
\frac14\nu_1\eta\varphi_1,\quad \mu\be(g(u)-g(w)-g'(w)\xi)\varphi_1\leq
\frac14\nu_1\xi\psi_1. \eeq
Combining \eqref{226}-\eqref{228} with again  the Poincar$\acute{e}$
inequality yields
\beq
\label{a15}
\frac{\mathrm{d}}{\dt}\int_{\Omega}\bigg(\frac{1}{2}\xi^2+\frac{1}{2}\eta^2+\big[\psi_1\xi+\varphi_1\eta\big]\bigg)\dx+\lam_1\|\xi\|_2^2
+\lam_1\|\eta\|_2^2+\frac{\nu_1}{2}\int_{\Omega}\big[\psi_1\xi+\varphi_1\eta\big]\dx
\leq 0.
\eeq
Let
$Y=\dis\int_{\Omega}(\xi^2+\eta^2+2\big[\psi_1\xi+\varphi_1\eta\big])\dx.$
Note that
\beq
\label{aaa8}
\int_{\Omega}\big[\psi_1\xi+\varphi_1\eta\big]\dx\geq 0. \eeq
By \eqref{a15} there holds \beq
\frac{\mathrm{d}Y}{dt}+\gamma Y\leq 0, \quad
\gamma=\min\bigg\{2\lambda_1,\frac{\nu_1}{2}\bigg\}, \eeq which again
yields $Y\leq Y(0)e^{-\gamma t}.$ Then by noting \eqref{aaa8} again
it follows
 that \beq\label{aaa9}
  \|u(x,t)-w(x)\|_2^2+\|v(x,t)-z(x)\|_2^2\leq Y(t)
  \leq C_0 \exp\bigg({-\min\bigg\{\lambda_1,\frac{\nu_1}{2}\bigg\}t}\bigg),
\eeq
for large $t$ with $Y(0)\leq C_0=\|w_1-w_{\lam,\mu}\|_2^2+\|z_1-z_{\lam,\mu}\|_2^2
+2\|\psi_1(w_1-w_{\lam,\mu})+\varphi_1(z_1-z_{\lam,\mu})\|_1.$\\

\textbf{Step 3: Consider more general $(u_{0},v_{0})$.}

Note the  comparison principle. If $(0,0)\le(u_0,v_0)\le(w_{\lam,\mu},z_{\lam,\mu}),$ then \eqref{215} holds. If $(w_{\lam,\mu},z_{\lam,\mu})\le(u_0,v_0)\le\not\equiv(w_1,z_1)$,  \eqref{aaa9} holds. At last, for more general $(0,0)\le(u_0,v_0)\le\not\equiv(w_1,z_1),$  \eqref{rate} holds for $C_0=\|w_{\lam,\mu}\|_2^2+\|z_{\lam,\mu}\|_2^2
+2\|\psi_1w_{\lam,\mu}+\varphi_1z_{\lam,\mu}\|_1+\|w_1-w_{\lam,\mu}\|_2^2+\|z_1-z_{\lam,\mu}\|_2^2
+2\|\psi_1(w_1-w_{\lam,\mu})+\varphi_1(z_1-z_{\lam,\mu})\|_1.$ The proof of Theorem \ref{Convergence rate} is therefore completed. \proofend

\appendix

\renewcommand{\appendixname}{Appendix}

\begin{appendices}
\section{}

\label{Stationary problem} We will show the proof of Theorem
$\textbf{A}$ in this Appendix. First we will prove that the elliptic
problem \eqref{E} has a classical solution for $\lam$ and $\mu$
small enough, while \eqref{E} has no solution for $\lam$ or $\mu$
large enough. More precisely, we will prove that the set \beq
\Lam:=\{(\lam,\mu)\in\r^+\times\r^+:\eqref{E} ~\mbox{has a classical
minimal solution}\} \eeq is nonempty and bounded.

\begin{lem}
\label{Lam bounded} $\Lam$ is bounded, and there exist $\lam_0>0$,
$\mu_0>0$ such that $(0,\lam_0]\times(0,\mu_0]\subseteq\Lam$.
\end{lem}

\proof Let $\gam\in H^1_0(\Om)$ be the regular solution of
$-\D\gam=1$ in $\Om$. It is then easy to verify that there exists
$\alp\in(0,{1}/{\|\gam\|_{\infty}})$ such that
$(\alp\gam,\alp\gam)$ is a supersolution of \eqref{E} if
$$\lam<\dis\f{1}{\sup\limits_{x\in\Om}\alp(x)}\sup\limits_{0<s<\f{1}{\|\gam\|_{\infty}}}
 \f{s}{f(s\|\gam\|_{\infty})}=:\lam_0\  \text{and}\
\mu<\dis\f{1}{\sup\limits_{x\in\Om}\beta(x)}\sup\limits_{0<s<\f{a}{\|\gam\|_{\infty}}}
\f{s}{g(s\|\gam\|_{\infty})}=:\mu_0.$$
As (0,0) is a subsolution and $\alp\gam>0$ in $\Om$, \eqref{E}
admits a regular solution for $\lam\in(0,\lam_0]$ and
$\mu\in(0,\mu_0]$. In fact, for these $\lam,\mu$, using
\eqref{assumption on initial data} and the monotone iteration for
$n\in\mathbb{N}$, \beq \left\{\arraycolsep=1.5pt
\begin {array}{lll}
w_0=z_0=0,&\\
-\Delta w_{n+1}=\lambda \al(x)f(z_n),\   \ &\hbox{in} \,\ \Omega,\\[2mm]
-\Delta z_{n+1}=\mu \be(x)g(w_n), \   \ &\hbox{in} \,\ \Omega,\\[2mm]
w_{n+1}=z_{n+1}=0,\ &\hbox{ on} \ \partial \Omega,
\end{array}
\right. \eeq we get the minimal solution
$(w_{\lam,\mu},z_{\lam,\mu})=\lim\limits_{n\rightarrow+\infty}(w_n,z_n)$.
Therefore, $\Lam$ is nonempty.

On the other hand, take a positive first eigenfunction $\varphi$ of
$-\D$ in $H^1_0(\Om)$ with the first eigenvalue $\lam_1$ such that
$\Io\varphi\dx=1$. By \eqref{E} and $w<1$, $z<1$, we arrive at \beq
\left\{ \bea{l} \lam_1\ge\lam_1\Io w\varphi\dx=\Io\varphi(-\D w)\dx
=\lam\Io \alp(x)f(z)\varphi\dx\ge\lam\Io \alp(x)f(0)\varphi\dx,\\
\lam_1\ge\lam_1\Io z\varphi\dx=\Io\varphi(-\D z)\dx =\lam\Io
\beta(x)g(w)\varphi\dx\ge\mu\Io \beta(x)g(0)\varphi\dx. \eea \right.
\eeq So $\Lam$ is bounded and
$\Lam\subseteq\big(0,\f{\lam_1}{\io
\alp(x)f(0)\varphi\dx}\big]\times\big(0,\f{\lam_1}{\io
\beta(x)g(0)\varphi\dx}\big]$. \proofend

Denote $\mu=\Gamma(\lam)$ as the critical curve such that if
$0\leq\mu<\Gamma(\lam),$ then $(\lam,\mu)\in\Lam;$ if
$\mu>\Gamma(\lam),$ then $(\lam,\mu)\in
(\mathbb{R}^+\times\mathbb{R}^+)\setminus \bar{\Lam}.$ By Lemma
\ref{Lam bounded}, there further hold  $0<\mu^*:=\Gamma(0)<+\infty$
and $0<\lam^*:=\Gamma^{-1}(0)<+\infty.$

Next we state that the critical curve $\mu=\Gamma(\lam)$ is
non-increasing. More precisely,
\begin{lem}\label{19}
If $0\leq\lam'\leq\lam, 0\leq\mu'\leq \mu$ for some
$(\lam,\mu)\in\Lam,$ then $(\lam',\mu')\in\Lam$.
\end{lem}
\proof Indeed, the solution associated to $(\lam,\mu)$ turns out to
be a super-solution to $(E)$ with $(\lam',\mu').$ \proofend
\\
\proofthA Define $\O_1=\Lam\backslash\Gam$. For
$(\lam_1,\mu_1),(\lam_2,\mu_2)\in\O_1$, there exist
$\theta_1,\theta_2>0$ such that $\mu_1=\theta_1\lam_1$ and
$\mu_2=\theta_2\lam_2$. Using Lemma \ref{19}, we can define a path
linking $(\lam_1,\mu_1)$ to $(0, 0)$ and another path linking
$(0,0)$ to $(\lam_2,\mu_2)$, which implies that $\O_1$ is connected.
Now, define $\O_2=(\R^+\times\R^+)\backslash\{\Lam\bigcup\Gam\}$.
Let $(\lam_1,\mu_1),(\lam_2,\mu_2)\in\O_2$. Then by Lemma \ref{19} again that $(\lam_{\max},\mu_{\max})\in\O_2$, where
$\lam_{\max}=\max\{\lam_1,\lam_2\}$ and
$\mu_{\max}=\max\{\mu_1,\mu_2\}$. We can take a path linking
$(\lam_1,\mu_1)$ to $(\lam_{\max},\mu_{\max})$ and another path
linking $(\lam_{\max},\mu_{\max})$ to $(\lam_2,\mu_2)$, which
follows that $\O_2$ is connected.

At last, it is reduced to prove that problem \eqref{E} admits no
weak solution for $(\lam,\mu)\in \O_2$. Suppose on the contrary that
$(w,z)$ is a weak solution to \eqref{E}. By the monotonicity of
$f,g$, it is easy to verify that for any $\del>1$,
$(\hat{w},\hat{z})=({w}/{\del},{z}/{\del})$ is a weak
super-solution for problem
$$\left\{\arraycolsep=1.5pt
\begin {array}{lll}
-\Delta w=\displaystyle\frac{\lambda}{\del}\alpha(x)f(z),\ \ &\hbox{in} \,\  \Omega,\\[2mm]
-\Delta z=\displaystyle\frac{\mu}{\del}\beta(x)g(w), \ \ &\hbox{in} \,\  \Omega,\\[2mm]
w=z=0,\ &\hbox{on} \ \partial \Omega,
\end{array}\right.\eqno{(E_{\del})}
$$
then the monotone iteration will enable us a weak solution
$(\tilde{w},\tilde{z})$ of $(E_{\del})$ satisfying
$0\le\tilde{w}\le\hat{w}\le{1}/{\del}<1$, and
$0\le\tilde{z}\le\hat{z}\le{1}/{\del}<1$. The regularity theory
implies that $(\tilde{w},\tilde{z})$ is a regular solution of
$(E_{\del})$. This means that
$({\lam}/{\del},{\mu}/{\del})\in\O_1\bigcup\Gam$. Let $\del$ tend
to 1, we get $(\lam,\mu)\in\O_1\bigcup\Gam$, which  contradicts with
the assumption. Therefore, no weak solution exists for
$(\lam,\mu)\in\O_2$
and  the proof of Theorem $\bf{A}$ is completed. \proofend

\section{}
\label{proof of local existecnce and uniqueness}

In this Appendix, we will show the proof of Theorem B.

\proofthB We first show the uniqueness of the solution to \eqref{P}.
For any given $0<T_0<T,$ suppose $(\tilde{u},\tilde{v})$,
$(\hat{u},\hat{v})$ are two pair classical solutions of \eqref{P} on
the interval $[0,T_0]$ such that
$\|\tilde{u}\|_{L^{\infty}(\Omega\times[0,T_0])}<1,$
$\|\tilde{v}\|_{L^{\infty}(\Omega\times[0,T_0])}<1,$
$\|\hat{u}\|_{L^{\infty}(\Omega\times[0,T_0])}<1,$
$\|\hat{v}\|_{L^{\infty}(\Omega\times[0,T_0])}<1$.

Indeed, the difference $(U,V)=(\tilde{u}-\hat{u},\tilde{v}-\hat{v})$
satisfies
\begin{equation}
\left\{\arraycolsep=1.5pt
\begin {array}{lll}
U_t-\Delta U=\lam\al(x)f'(\th_v)V,\   \ &\hbox{ in} \,\  \Omega\times(0,T_0],\\[2mm]
V_t-\Delta V=\mu \beta(x)g'(\th_u)U,\   \ &\hbox{ in} \,\  \Omega\times(0,T_0],\\[2mm]
U=V=0,\ &\hbox{ on} \ \partial \Omega\times(0,T_0],\\[2mm]
U(x,0)=V(x,0)=0, \ &\hbox{ for}\ x\in\bar{\Omega}
\end{array}\right.
\end{equation}
where $\th_v$ is between $\tilde{v}$ and $\hat{v}$, $\th_u$ is
between $\tilde{u}$ and $\hat{u}$. The assumption on
$(\tilde{u},\tilde{v})$, $(\hat{u},\hat{v})$ implies that
$f'(\th_v),g'(\th_u)\in L^{\infty}(\Omega\times [0,T_0])$ for any
$T_0<T$. By using the comparison principle stated in Lemma \ref{MP},
we deduce that $U=V\equiv0$ on $\bar{\Omega}\times [0,T_0]$.

To obtain Theorem B, it is reduced to show the existence. Let
$(\zeta,\rho)$ be the solution of the  ODE system
\beq
\label{bbbb1} \left\{\arraycolsep=1.5pt
\begin {array}{lll}
\F{\mathrm{d}\zeta}{\dt}=\lam\|\al\|_{\infty}f(\rho),\   \ &\hbox{in} \,\  (0,T),\\[2mm]
\F{\mathrm{d}\rho}{\dt}=\mu\|\be\|_{\infty}g(\zeta), \   \ &\hbox{in} \,\  (0,T),\\[2mm]
\zeta(0)=\|u_0\|_{\infty},\ \rho(0)=\|v_{0}\|_{\infty}.&
\end{array}
\right. \eeq
The
local existence of \eqref{bbbb1} can be obtained by \cite[Chapter
III]{Wa}. Obviously, $(\zeta,\rho)$ is a supersolution of \eqref{P}.
Since $(0,0)$ is a subsolution of \eqref{P}, it follows from
\cite[Theorem 2.2]{Pao2} that there exists a unique classical
solution $(u,v)$ to \eqref{P} between $(0,0)$ and $(\zeta,\rho)$. In
conclusion, the proof of Theorem B is completed. \proofend

 \end{appendices}

 \section*{Acknowledgements}
The authors sincerely thank Professor Rodrigo Clemente, Professor Marcelo Montenegro and Professor Dong Ye for their kind help.

\end{document}